\documentclass{article}
\usepackage[greek,french,english]{babel}
\usepackage{csquotes}
\usepackage{graphicx}
\usepackage{csquotes}
\usepackage{dsfont}
\usepackage{newtxtext}
\usepackage{newtxmath}
\begin{document}

\title{The Tautochrone of Huygens and Abel:\\ From Constructive Geometry 
\\ to Fractional Calculus}
\author{Luiz Roberto Evangelista$^{1,2,3}$ \and Francesco Mainardi$^4$}
\date{$^1$Departamento de F\'isica, Universidade Estadual de Maring\'a,  Maring\'a,  Paran\'a, 87020-900, Brazil.
\\ $^2$Istituto dei Sistemi Complessi del Consiglio Nazionale delle Ricerche (ISC--CNR), Via dei Taurini 19,  00185  Roma, Italy.
\\ $^3$Department of Molecular Science and Nanosystems, 
Ca' Foscari University of Venice, Via Torino 155,  30175 Mestre (VE),  Italy. 
\\
$^4$Department of Physics and Astronomy, University of Bologna, and INFN,
Via Irnerio 46, 40126 Bologna, Italy. \\[2ex]%
\today\\
Eur. Phys. J. Plus (2026) 141:557\\
https://doi.org/10.1140/epjp/s13360-026-07786-3\\
Published on line, 19 May 2026 
}

\maketitle


\begin{abstract}
In this paper, we explore the connections between Christiaan Huygens and Niels Henrik Abel through the tautochrone problem. The problem --- determining the curve along which a particle descends under gravity in the same time, regardless of its starting point --- has been a central topic at the intersection of physics, geometry, and analysis. Though these two major figures are separated by nearly two centuries, they approached the problem in radically different ways. While Huygens proposed a physical solution based on geometric construction, Abel approached the problem within the analytic framework of integral equations, employing a procedure that can be seen as  anticipating and paving the way for the development of differential calculus of arbitrary order. 
This contrast highlights a broader historical narrative: the transformation of mathematical thinking from constructive geometry to abstract analysis.
\end{abstract}

\noindent Keywords: Isoperimetric problems; variational calculus; integral equations; fractional calculus

\section{Introduction}
\label{Sec:Introd}

Problems involving extremes (maximums and minimums) are old and persistent. The very use of terms like ``law of least effort'', ``path of least resistance'', ``shortest distance between two points'', etc., is part of our common sense of economy or simplicity in the actions of nature~\cite{Fermat}. The first problem of finding an extreme that is known to us is associated with the name of the Phoenician Princess Dido, queen of Carthage --- immortalized by the Latin poet Virgil in the \textit{Aeneid}, and later by the English composer Henry Purcell (1659--1695) in the musical drama \textit{Dido \& Aeneas}. The legend says that the princess, having emigrated to North Africa after being defrauded of her belongings by her brother, Pygmalion, was promised by Hiarbas, the native chief of the region, as much land as she could enclose in a cowhide. Cleverly, she cut the cowhide into a very long strip and placed it between herself and the sea, creating a territory large enough to build the city of Carthage~\cite{Kelvin}. The typical solution found by Dido is to arrange the strip in a circular shape. It is, therefore, an \textit{isoperimetric problem}, one of many others that will be encountered throughout history. Another problem of this nature is associated with the name of Heron of Alexandria (10--70~A.D.). He is credited with deriving the law of reflection from the principle that the light ray emitted from point A reaches point B after being reflected by a mirror, following the path that takes the least time. This principle was extended by Fermat to derive Snell--Descartes' law of refraction~\cite{Fermat}.

In this work, we focus on a problem of this nature, which also belongs to the realm of isoperimetric and optimization problems --- the \textit{tautochrone problem}. This problem concerns the determination of the curve along which a particle, starting from any point on the curve under the gravitational force, will reach the lowest point in the same amount of time, regardless of its initial position. The tautochrone problem is closely related to the concept of extremal times and optimal paths, much like the problems that Dido and Heron of Alexandria faced. In fact, solving the tautochrone problem is akin to finding an optimal path that minimizes a certain quantity --- in this case, the time of descent --- given constraints. It's another example of how nature often seeks optimization, whether through the geometry of curves or the behavior of physical systems under certain conditions.

 The problem was also posed by the mathematician Christiaan Huygens in the 17th century and was later shown to have a solution in the form of the \emph{cycloid}. What Huygens demonstrated was that a material point, released from rest and allowed to slide without friction along an arc of an inverted cycloid, reaches the lowest point of the curve in the same amount of time, regardless of the initial height. Thus, a pendulum constrained to oscillate along such a cycloidal arc exhibits a period of oscillation that is truly independent of the amplitude of its motion. Huygens' approach was fundamentally constructive. He did not yet have access to the fully developed calculus of Newton and Leibniz, but relied instead on geometric arguments and proportions, integrating physical intuition with mathematical ingenuity. His solution was elegant, self-contained, and tailored to a particular problem with practical consequences.

By contrast, in the early 19th century, Niels Henrik Abel approached a more general class of problems: inverse integral equations. In a brief but brilliant life, Abel transformed how mathematicians viewed such problems. He considered equations of the form
\begin{equation}
\psi(a) = \int_{x=0}^{x=a} \frac{f(x)}{\sqrt{a - x}} \, dx,
\end{equation}
and posed the question: given $\psi(a)$, can one determine the function $f(x)$? This type of integral, now known as an Abel integral, has direct analogies to the time-of-descent problem in the tautochrone curve, where the time taken to reach a point is an integral over a function of the height difference. Abel's treatment of this inverse problem marked a turning point. Instead of seeking a specific curve, he found a general method for recovering unknown functions from their integrals --- paving the way for later developments in functional analysis and integral transforms. The tautochrone problem became, in his hands, a special case of a much broader theory.

Here, we aim to highlight the scientific contributions of these two mathematicians, using the tautochrone problem as a guiding thread. Some biographical details are included to underscore the contrast between their backgrounds, beginning in childhood, when both had already shown a rare and precocious talent for mathematical science and abstract thinking. By focusing on how each of them approached the problem --- separated by two centuries --- we also trace how analysis became increasingly sophisticated in addressing important and advanced problems in physics and mathematics. Special attention is given, in the final part of the paper, to the manners in which Abel's work paved the way for the development of key concepts in fractional calculus --- matured in the years that followed~\cite{Ross1977}.

\section{The Geometric Construction of Huygens}
\label{Sec:Huygens}

Christiaan Huygens was born on April 14, 1629, in \textit{Den Haag} (The Hague), in the Dutch Republic. He was the son of Constantijn Huygens, a diplomat with a strong background in philosophy and the natural sciences. Christiaan's father was also a distinguished poet, securing an enduring place in the history of Dutch literature. Constantijn ensured that his son received an exceptional education and had access to the leading intellectual circles of the time. A prominent role in this upbringing is attributed to Father Marin Mersenne (1588--1648) --- who corresponded regularly with Huygens's father --- as well as to the friendship the elder Huygens maintained with Ren\'e Descartes (1596--1650). Educated at home by private tutors until the age of sixteen, Christiaan benefited from the intellectual environment surrounding his family. Descartes, who was then living in the Dutch Republic and maintained a friendly relationship with Constantijn, is believed to have taken interest in the young man's mathematical development. In 1645, Huygens enrolled at the University of Leiden, where he studied law, mathematics, and classical languages. Two years later, in 1647, he continued his legal studies at the College of Orange in Breda. The period from 1650 to 1666 would prove to be the most fertile and productive phase of his scientific life.

\subsection{The Tautochronism of the Cycloid}

Because astronomical investigations require increasingly accurate time measurements, Huygens devoted considerable effort to the problem. In 1656, he patented the first pendulum clock, significantly improving the precision of timekeeping. This innovation was later described in his 1658 publication \textit{Horologium} (``The Clock''). Although his investigations into centrifugal force also date back to this period (1659), it was in his major work of 1673, \textit{Horologium Oscillatorium sive de motu pendulorum} (``The Pendulum Clock, or On the Motion of Pendulums''), that he developed a full theoretical account of pendular motion. It is also in this treatise that the modern formula for the centrifugal force in uniform circular motion first appears. The \textit{Horologium Oscillatorium} is often regarded as one of the most important work on mechanics prior to Newton's \textit{Principia} (1687). In the book, Huygens derived the formula for the period of a simple pendulum under the small-angle approximation --- a classical result encountered in the Galileo's work. Huygens went further by investigating the precise curve along which a pendulum would swing with a period that is truly independent of its amplitude. His solution identified this path as a cycloid.

This remarkable property is known as the \textit{tautochronism} of the cycloid --- from the Greek words \textgreek{taut\'o} (equal, same) and \textgreek{kr\'onos} (time) ---  and it requires that the pendulum follow an inverted cycloidal trajectory. The tautochrone problem is addressed in Part II, specifically in Propositions XVI through XXVI. Let us now proceed by sketching some of the representative steps of his derivation, focusing the last two propositions, starting with Proposition XXV. It states that~\cite{Huygens1673}:
\begin{center}
\textsc{Proposition XXV}
\end{center}
{\textit{In a cycloid with a vertical axis, and with the vertex seen to be the lowest point,
the times of descent for some body, on leaving any point on the cycloid from rest until it
reaches the lowest point at the vertex, are equal to each other; and this time has the
same ratio to the time of fall along the whole axis of the cycloid as the semi-circumference of the circle to the diameter.}}

\begin{figure}
\centering
\includegraphics[scale=.60,angle=0]{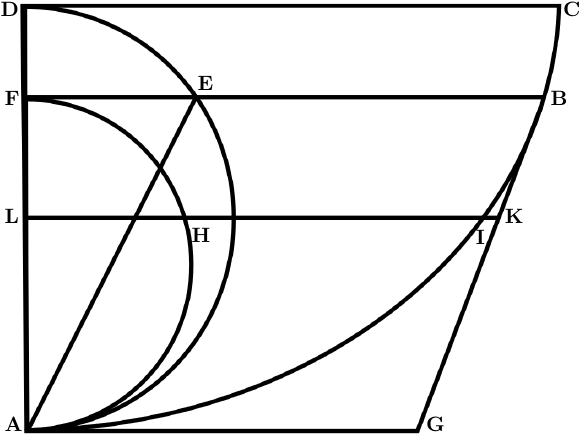}
\caption{This is the original figure presented in the Huygens work in connection with the results analyzed in Propositions XXV--XXVI. The particle is released from the point B and the final point of its trajectory is in A. We notice that EA is parallel to BG and both describe a kind of inclined plane for the descend of the particle. }\label{Huygens-25}
\end{figure}

To reproduce his construction, we need the help of Fig.~\ref{Huygens-25}, part of the Huygens' book~\cite{Huygens1673b}.
For a better understanding of the discussion that follows, here is a summary of the elements that appear in Fig.~\ref{Huygens-25} and their meaning:
\begin{itemize}
\item Point A --- is the lower vertex of the cycloid, the endpoint of the motion;
\item BA is an arc of the cycloid --- the path followed by a body in free fall;
\item Vertical axis DA --- a straight line along which the `standard' time is measured for comparison of the motion;
\item Semicircle FHA --- constructed on the diameter FA, used to establish the proportion of times;
\item Tangent BG --- the line tangent to the cycloid at point B, parallel to the line EA, which intersects the semicircle.
\end{itemize}

The cycloid shall be ABC the vertex of which is seen to be A pointing downwards,
with the vertical axis truly AD, and by taking some point on the cycloid, such as B, a
body is free to fall from natural forces along the arc BA, or along a surface thus curved
into this shape. Huygens says that the time of descent of this body to be as to the time of fall along
the axis DA, thus as the semi-circumference of the circle to the diameter. From which
demonstration, the time of descent along any arc of the cycloid to finish at A also, are
agreed to be equal to each other~\cite{Bruce}. Huygens then states that the descent time of this body is, in relation to the fall time along the DA axis, as the semicircle is in relation to the diameter. This proportion stated by Huygens --- between the descent time on the cycloid and the vertical fall time along the axis --- is the expression of tautochronism: in fact, the proportionality constant is $\pi/2$, since the perimeter of the semicircle is $\pi\,r$, while its diameter is $2r$.

The discussion continues by diving into the details of the geometric construction. Huygens explains that a semicircle is drawn on the axis DA. This semicircle intersects the line BF at point E, which is parallel to the base DC. The line EA is then drawn, and the line BG is drawn parallel to it, which turns out to be a tangent to the cycloid at point B. In fact, this line also intersects the line drawn horizontally from point A at point G. Additionally, a semicircle FHA is drawn on the line FA. With this setup, Huygens states that the time it takes to descend along the arc of the cycloid BA is to the time of uniform motion along the line BG (with half the speed of 
BG) as the arc of the semicircle FHA is to the straight line FA.  In other words, the time for uniform motion along BG is the same as the time it would take to descend along the same line BG, or along EA (which is parallel and equal to BG), or equivalently, the time to fall under gravity along the axis DA,  as shown in Proposition 6 by Galileo~\cite{Galileo}. Therefore, the time along the arc BA is related to the time of descent along the axis DA in the same way the circumference of the semicircle FHA is related to the diameter FA.

The last proposition dedicated to this problem also refers to Fig.~\ref{Huygens-25} and is stated as follows~\cite{Bruce-2}:
\begin{displayquote}
\begin{center}
\textsc{Proposition XXVI}
\end{center}
\textit{With the same positions, if some line HI is drawn above which cuts the arc BA in I,
and the circumference FHA in H: I say that the time to pass through the arc BI, to the
time to cross the arc IA after BI, has the ratio which the arc of the circumference FH
has to HA.}
\end{displayquote}

The Huygens' explanation of the construction may be paraphrased as follows~\cite{Huygens-1673-b}.  Consider the straight line HI intersecting the tangent BG at point K, and the axis DA at point  L. The time it takes to travel along the arc BA is compared to the time it takes for uniform motion along BG with half the speed of BG. This comparison is similar to how the arc FHA compares to the straight line FA.  The time of uniform motion along BG is compared to the time of uniform motion along BK, with the same half-speed of BG. The ratio of these times is the same as the ratio of BG to BK, which is equivalent to the ratio of FA to FL.  The time for uniform motion with the same speed along BK is compared to the time it takes to move along the arc BI, which is similar to how the length FL compares to the arc FH. 

As a result, the time to travel along the arc BA is equal to the time to travel along the arc BI, and this is in the same proportion as the arc FHA is to the arc FH. Finally, by dividing and converting, the ratio of the time along BI to the time along IA after BI is the same as the ratio of the arc FH to HA.  We underline that both propositions end by the Q.E.D. statement in the original work because for Huygens the geometrical constructions he employs are rigorous method to demonstrate the validity of a physical result. This is indeed the case as we can show hereafter by recalling the problem using simpler language and more updated mathematical notation.

To proceed this way, let's first remember that Huygens started with the geometric generation of the cycloid, i.e., he considered a circle of radius $r$ rolling on a horizontal surface without slipping. A fixed point on the edge of the circle traces a path that defines the cycloid. The resulting curve is symmetric, with its lowest point at the vertex. When the curve is inverted (with the concavity facing upwards), it gives the ideal profile of a tautochrone pendulum. Huygens' analysis leads us to think that, for the descent time along a curve to be independent of the starting point, the curve must have the property of equating the acceleration of a body along it with the appropriate free fall acceleration for comparison purposes.

\subsection{The Tautochrone as a Variational Problem}

Using a modern approach, we could say that he wants to find a curve such that
\begin{equation}
\int \frac{ds}{\sqrt{2gy}} = \text{constant},
\end{equation}
where $ds$ is an element of the curve, $g$ is the acceleration due to gravity, and $y$ denotes the vertical axis in a Cartesian system. In this approach, the deduction is based on comparing two motions, namely, the motion along the inverted cycloid (frictionless sliding) and uniform circular motion, i.e., motion with constant speed along a circle. Based on these premises, Huygens establishes a correspondence between the vertical height of the particle on the cycloid and the angular position of a point moving uniformly along a circle. He then shows that, if one considers an auxiliary circle that generates the cycloid (by rolling), the time it takes for a point to reach a certain height on the cycloid coincides with the time a point in uniform motion takes to reach the same angle on the circle --- provided the variables are correctly associated.
 In a few words, the projection of uniform circular motion along the generatrix of the cycloid coincides with the free fall motion along the cycloid --- and since uniform circular motion is isochronous, the motion along the cycloid will also be isochronous!

In proceeding, Huygens notes that a common pendulum (one suspended by an ideal string) is not isochronous: its period depends on the amplitude of the oscillation. To resolve this, he proposes restricting the motion of the string so that the pendulum, when oscillating, follows an inverted cycloid rather than a circular arc. This ensures that, regardless of the amplitude, the oscillation period of the pendulum remains exactly constant. With these ingredients, Huygens then demonstrates that the inverted cycloid is the only curve for which the descent time is independent of the starting point. In fact, the descent time from any point to the vertex is
\begin{equation}
T = \pi \sqrt{\frac{r}{g}},
\end{equation}
which coincides with the period of a cycloidal pendulum. Finally, the ratio between this time and the free fall time along the cycloid's axis is given by

\begin{equation}
\frac{T_{\text{cycloid}}}{T_{\text{free \, fall}}} = \frac{\pi}{2}.
\end{equation}

In other words and somewhat more formal terms, we wish to find the curve $y(x)$ such that
\begin{equation}
\int_A^B \frac{ds}{v(y)} = \text{constant},
\end{equation}
where $ v(y) $ is the velocity of the particle along $ y $. In this case, we have
\begin{equation}
ds = \sqrt{1 + \left(\frac{dy}{dx}\right)^2}dx = \sqrt{1 + y'^2}dx
\end{equation}
and, also, the velocity is given by $ v(y) = \sqrt{2gy} $, a result derived using the conservation of the mechanical energy. Thus, the descent time for the particle can be defined as
\begin{equation}
T = \int \frac{\sqrt{1 + y'^2}}{\sqrt{2gy}} dx,
\end{equation}
which is the same integral that appears in the \textit{brachistochrone} problem --- the curve of minimum time between two points. 

The starting point of this investigation, in fact, was a challenge posed by Johann Bernoulli in 1696. It consisted of finding the trajectory followed by a body starting from rest at a given point, in a gravitational field, along a frictionless curve, until the final point. The curve is known as the brachistochrone, from Greek: \textgreek{br\'akistos}     (shortest) and  \textgreek{kr\'onos} (time). The problem was tackled (and solved) by Jakob I Bernoulli, in an article dedicated to solving one of the problems proposed by his brother: \textit{Solutio problematum fraternorum} (Solution to the brother's problem), which he titled \textit{De curva celerrimi descensus} (literally: the curve of fastest descent)~\cite{Bernoulli}. Euler later solved the same problem, thereby founding the calculus of variations~\cite{Euler1744,Fraser1994}. Later, Lagrange revolutionized the field by introducing the $\delta$ algorithm and incorporated variational methods into a much more general mathematical framework~\cite{Lagrange-LF,Fraser-LF}.

As can be seen from the excerpts of Huygens' demonstration, none of these equations were employed in establishing his results. We have introduced them here merely to bring Huygens' language closer to the one we use in the study of the fundamental problems of mechanics. This fact, by itself, illustrates the importance that analytical methods have played in mechanics. Indeed, in the period between the publication of Newton's \textit{Principia} in 1687 and the emergence of Lagrange's \textit{M\'echanique Analytique} in 1788, a highly elaborate language in mathematical physics developed --- an outcome of the effort to integrate mechanics into mathematical analysis. It was a fertile period during which scientists of great stature, such as Leibniz, Varignon, Jacob I and Johann I Bernoulli (and his son Daniel), Euler, d'Alembert, Laplace and Lagrange, among others, devoted themselves wholeheartedly to this endeavor. The work of these men led to the emergence of analytical mechanics, practically in the form we still know and that is presented in our textbooks today.

\section{The Integral Equations of Abel}

Niels Henrik Abel was a Norwegian mathematician renowned for his groundbreaking work in algebra. He is best known for proving that the general equation of the fifth degree (quintic equations) cannot be solved by radicals, a result now known as Abel's impossibility theorem (Abel--Ruffini theorem)~\cite{Abel1824,Abel1826b}. Abel's work laid the foundation for much of modern algebra and group theory. Despite dying young at 26, Abel made significant contributions to mathematics, including studies on elliptic functions and the development of Abelian integrals. His legacy is celebrated in the mathematical community, and the prestigious Abel Prize is named in his honor.

Abel was born on August 5, 1802, in Nedstrand, Norway, into a modest family. His father, a Lutheran pastor, struggled with alcoholism, which led to financial difficulties. Despite this, Abel showed early signs of exceptional intelligence. He received his education at the local school in his hometown before moving to the University of Christiania (now Oslo) in 1821. Abel was initially enrolled in theology, as his family hoped he would follow in his father's footsteps as a priest. However, he soon turned to mathematics, developing a deep interest in the subject, particularly in the areas of algebra and functions. His talent in mathematics became clear to his professors, and he quickly gained a reputation as a prodigy. Abel's time at the university was financially strained, and he often had to struggle with limited resources. Despite this, he completed his studies, and by 1823, he sent his first important mathematical paper to the renowned mathematician Carl Friedrich Gauss. Gauss was impressed by Abel's work, leading to a correspondence that would later have a significant influence on his career. Abel's education was largely self-directed, supplemented by his interactions with prominent mathematicians of the time. However, due to financial hardship, he was unable to secure a formal academic position. As a result, he spent the last years of his life in poverty, moving between various cities in Europe, where he worked on his mathematical research. In 1825, Abel moved to Paris, where he met several important mathematicians, including Joseph Fourier, and continued his studies. Unfortunately, his health began to deteriorate, and he died of tuberculosis on April 6, 1829, at the age of 26. Despite his short life, Abel's contributions to mathematics were profound and continue to influence the field to this day. His early life was marked by personal and financial struggles, but his mathematical brilliance shone through, leading to the development of key theories in algebra and analysis. This is completely different from the early life of Christiaan Huygens, who grew up in good health, in a financially comfortable and intellectually rich household, where the atmosphere was stable, cultured, and conducive to learning.

In 1823, Abel solved the tautochrone problem by using a fractional integration of order $\alpha = 1/2$~\cite{Abel1823} and, in 1826, extended the approach to $\alpha \in (0, 1)$~\cite{Abel1826}.
In his 1826 work, Abel presents a different solution method, based on the use of the properties of Euler's gamma function. As emphasized in Ref.\cite{Podlubny-Magin}, it is a more concise and elegant method which, however, abandons the path opened toward fractional calculus, characteristic of his 1823 paper. Indeed, as we briefly discuss now, the solution proposed by him to the tautochrone problem involves the pioneering use of an integral operator of arbitrary order, and is the first explicit solution known of an integral equation\cite{Lutzen1982}.

\subsection{Abel's Solution: Remarkable Theorem}

On the first page of Abel's 1823 article, dedicated to definite integrals, we read~\cite{Abel1823}:

\begin{quote}
{\emph{It is well known that, with the help of definite integrals, many problems can be solved which otherwise cannot be resolved, or at least are very difficult to handle [...] I will demonstrate a new application by solving the following problem.}}
\end{quote}

\begin{figure}[htp]
\centering
\includegraphics*[scale=.650,angle=0]{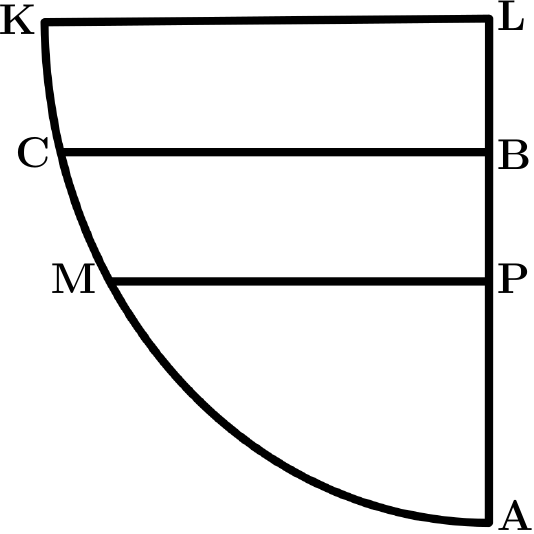}
\caption{The illustrative picture appearing in the French translation of the Abel's 1823 paper~\cite{Abel1823}. }
\label{Abel-1}
\end{figure}

Abel then considers a problem more general than the classical tautochrone, which we now briefly summarise in order to follow the core ideas and consequences of his approach. Consider Fig.~\ref{Abel-1}, which illustrates the motion of a particle along an arc, such as the one described by the curve CA. Let CB be a horizontal line and AB perpendicular to CB. The curve AM is described in Cartesian coordinates, where $AP = x$ and $PM = y$. We also define $AB = a$ and denote by $AM = s$ the arc length along the curve, which is initially unknown.

The particle begins its motion from rest at point C and moves along the curve CA. The time $T$ required to traverse the arc depends both on the shape of the curve and on the distance $a$. The problem is thus to determine the form of the curve KCA such that the time of descent $T$ is a given function of $a$, for instance, $T = \psi(a)$. We use today's notation; the  notation employed by Abel is such that a function $\psi(a)$ is denoted by  $\psi\,a$.

After some mechanical considerations, Abel found for this problem the following equation defining the time involved in the motion of the particle:
\begin{equation}
\psi(a) = \int_{x=0}^{x=a} \frac{ds}{\sqrt{a-x}},
\end{equation}
in which $s=s(x)$ is the unknown function. However, instead to handle with it directly, Abel considers the more general equation:
\begin{equation}
\label{Abel-eq-1}
\psi(a) = \int_{x=0}^{x=a} \frac{ds}{(a-x)^n},
\end{equation}
where $n$ has to be less than 1 to prevent the infinite integral between two limits; $\psi(a)$ is an arbitrary function that is not infinite, when $a=0$. 

To proceed,  let us recall that the problem is to find the curve $s(x)$ which Abel searches in the form of a power series
\begin{equation}
s (x) = \sum_m a^{(m)} x^m.
\end{equation}
After term-by-term operations with the series, and after using the properties of the gamma function --- by quoting in this regard the work of A. M. Legendre ---, and after some pages of manipulations with the power series, he arrives at the following solution to Eq.~(\ref{Abel-eq-1}), called by him as a ``remarkable theorem''~\cite{Podlubny-Magin}:
\begin{equation}
\label{Eq-remark}
s(x) = \frac{\sin n \pi}{\pi} x^n \int_0^1 \frac{\psi(xt) dt}{(1-t)^{1-n}}. 
\end{equation}
In a few words, the remarkable result is that \textit{given an integral like the one in Eq.~(\ref{Abel-eq-1}), the curve $s(x)$ is given by Eq.~(\ref{Eq-remark})}.  After examining the obtained formulas from other perspectives, Abel mentions that $s(x)$ can be expressed in a different way, namely:
\begin{equation}
s(x) = \frac{1}{\Gamma (1-n)}\int^n \psi(x) dx^n = \frac{1}{\Gamma (1-n)} \frac{d^{-n}}{dx^{-n}}\psi(x). 
\end{equation}
As underlined in Ref.~\cite{Podlubny-Magin}, it is easy to recognize, behind the archaic notation he used, that Abel presents two expressions for the fractional-order integral: a derivative of negative order and a symbol that later appears in Liouville's works for the integral of order $n$: $\int^n \psi(x) dx^n$ (See Refs.~\cite{Liouville1832b,Liouville1832a,Liouville1834}).

Subsequently, Abel considers the case $n=1/2$, in order to obtain from the integral
\begin{equation}
\psi(a) = \int_{x=0}^{x=a} \frac{ds}{\sqrt{a-x}},
\end{equation}
the desired curve $s(x)$, when the time is given by $psi(a)$, i.e., 
\begin{equation}
s(x) = \frac{1}{\sqrt{\pi}} \frac{d^{-1/2}}{dx^{-1/2}}\psi(x) = \frac{1}{\sqrt{\pi}}\int^{1/2} \psi(x) dx^{1/2}, 
\end{equation}
from which one obtains:
\begin{equation}
\psi(x) = \sqrt{\pi} \frac{d^{1/2}}{dx^{1/2}}s(x). 
\end{equation}
For our analysis, this is a crucial point in the paper because Abel inverts the fractional integral by using a fractional differentiation on both sides of the equation, thus assuming (and understanding) that fractional order integration and differentiation are mutually inverse procedures. He affirms indeed~\cite{Abel1823}:
\begin{quote}
{\emph{If the equation of a curve is $s = \psi(x)$, the time a body takes to travel along an arc of the curve whose height is $a$ is $\sqrt{\pi} \frac{d^{1/2} \psi(a) }{da^{1/2}}$.}}
\end{quote}

To conclude this incursion for the original approach of Abel, we remark for future purposes that 
Eq.~(\ref{Abel-eq-1}) could be rewritten as:
\begin{equation}
\label{preCaputo}
\psi(t) = \int_0^t \frac{s'(x) }{(t-x)^n}dx,
\end{equation}
which, apart from a constant coefficient, may be identified with the  well-known 
Caputo's time-fractional derivative of real order to be introduced later~\cite{Caputo1967,Caputo1969}. In this perspective, the solution proposed by Abel is the inverse operation of the fractional integral of the same arbitrary order $n$~\cite{Podlubny-Magin}.

\subsection{Abel's Isochrone: Explicit Solution}

\begin{figure}
\centering
\includegraphics[scale=0.75, angle=0]{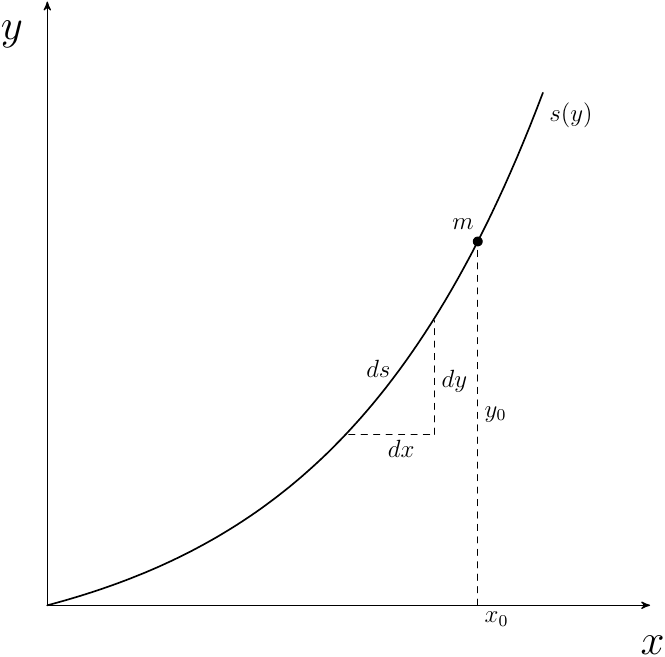}
\vspace{-0.25cm}
\caption{A mass $m$ under gravity slides down along a frictionless wire,  whose figure is defined by $s(y)$, such that $s(0)=0$. This trajectory is equivalent to the arc BA in Fig.~\ref{Huygens-25}, rebuild here with more appropriate notation.}
\label{Fig:Abel}
\end{figure}

To appreciate the result obtained by Abel in more simple terms, let us first reformulate the physical problem of a particle of mass $m$ sliding the curve shown in Fig.~\ref{Fig:Abel}. It is released from the rest at the position $y=y_0$ and slides under the influence of the gravity, along a frictionless wire.
The problem is to determine the curve $s(y)$, followed by the particle,  such that the time employed to slide down to the lowest point under gravity is independent of its initial position ($x_0, y_0$) on the curve. If the mass is released from $y=y_0 >0$, then its speed at $y$ will be given by
\begin{equation}
\label{Eq**Energia}
\frac{1}{2} m \left( \frac{d s}{dt} \right)^2 = m g y_0 - m g y = mg (y_0-y),
\end{equation}
from which we obtain:
\begin{equation*}
dt = \pm \frac{ds}{\sqrt{2g (y_0-y)}} \quad \longrightarrow \quad dt = - \frac{1}{\sqrt{2g (y_0-y)}}\, \frac{d s}{dy}\, dy,
\end{equation*}
because $ds/dy >0$, but $dy/dt <0$. The time of descent is then
\begin{equation}
\label{Tzero}
T(y_0) = \int_0^{y_0} \frac{1}{\sqrt{2g (y_0-y)}}\, \frac{d s}{dy}\, dy.
\end{equation}
We observe that~(\ref{Tzero}) was obtained by invoking only the conservation of the mechanical energy, Eq.~(\ref{Eq**Energia}). It is an integral equation to be solved in order to get $f(y) = ds/dy$. Let us change a little the notation ($ y_0 \to y, \, y \to z$) to search for a curve $s(z)$ such that
\begin{equation*}
\sqrt{2 g} T(y) = \int_0^y \frac{1}{\sqrt{y-z}} \frac{d s}{dz}\, dz = k (y) = \textrm{constant}.
\end{equation*}
The formal problem is then to find $f(z) = ds/dz$, given $k(y)$:
\begin{equation}
\label{Abelfirst}
k(y) = \int_0^y  \frac{f(z)}{\sqrt{y-z}} dz =\int_0^y f(z) K(y,z) dz,
\end{equation}
which may be recognized as a convolution between $f(z)$ and a kernel $K(y,z) = (y-z)^{-1/2}$. Equation~(\ref{Abelfirst}) is an integral equation for $f(z)$ and may be solved by taking the Laplace transform of its both sides. Remembering that the Laplace transform is defined as
\begin{equation}
\label{Def:LaplaceTransform}
{\mathcal{L}}\{ f(t); s\}:= F(s) = \int_0^{\infty} e^{-s t} f(t) dt, \quad \quad s\in \mathbb{C},
\end{equation}
where $\mathbb{C}$ stands for the set of complex numbers, after its computation, we obtain
\begin{equation}
\frac{k(y)}{s} = F(s) \frac{\sqrt{\pi}}{s^{1/2}}
\quad \Longleftrightarrow \quad  F(s) = \frac{c_1}{s^{1/2}},
\end{equation}
where $c_1 = k(y)/\sqrt{\pi}$ is a constant. 
The inverse transform of $F(s)$ is
\begin{equation}
f(z) = c z^{-1/2},
\end{equation}
where $c$ is another constant. Remember that $f(z) = ds/dz$ or $f(y) = ds/dy$. Thus,
\begin{equation}
f(y) = \frac{ds}{dy} = \left[1+ \left(\frac{dx}{dy}\right)^2\right]^{1/2} = c y^{-1/2}.
\end{equation}
By putting $y= c^2 \sin^2 (1/2\theta)$, a short calculation yields
\begin{equation}
x= \frac{1}{2}c^2 (\theta + \sin\theta) \quad {\textrm{and}} \quad
y = \frac{1}{2} c^2 (1-\cos\theta),
\end{equation}
which are the parametric equations of a cycloid~\cite{Sokolnikoff1958}. 

\subsection{The Birth of Fractional Calculus}

Let us pause for a moment to underline a central point of the present analysis. We remember that the integral equation obtained by Abel in the tautochrone problem can formally be written as
\begin{equation}
\label{Abel-A}
\int_0^x \frac{g(t)}{\sqrt{(x-t)^{1/2}}} dt = f(x)
\end{equation}
However, this is, essentially, the definition of a \textit{fractional integral of order $1/2$}, in the sense of Riemann--Liouvile, as we discuss soon.  We may rewrite it in the form:
\begin{equation}
{}_0{\textrm{I}}_x^{1/2} g(x) = \int_0^x \frac{g(t)}{\sqrt{(x-t)^{1/2}}} dt,
\end{equation}
where we have introduced the notation ${}_0{\textrm{I}}_x^{1/2}$ for an integration of arbitrary order $1/2$ in the variable $x$. Thus, the tautochrone problem in the Abel's perspective may assume the form:
\begin{equation}
f(x) = {}_0{\textrm{I}}_x^{1/2} g(x). 
\end{equation}
In this way as well, to obtain the unknown function $g(x)$ became possible if we promote the inversion of a fractional operator, i.e., of the fractional integral of order $1/2$.   In this framework, the Abel problem would be the particular case of the inversion of a fractional operator. The formal solution would require the application of the corresponding fractional derivative, i.e., 
\begin{equation}
g(x) = \frac{d^{1/2}}{dx^{1/2}}\int_0^x \frac{f(t)}{(x-t)^{1/2}} dt.
\end{equation}
This kind of simplified reasoning opens space to discuss the fractional calculus as a natural generalization of the classical differential operators, motivated by historical problems. Thus, the modern interpretation of the Abel problem as involving a fractional operator is an instructive bridge between the early years of the mechanics and mathematical analysis and the present. 

Let us now move to a more formal way to connect again these results from Abel to the roots of the fractional calculus.  To accomplish this task, the  problem can be reformulated in such a way that one has to find $f(z) = ds/dz$, given $k(y)=k$, such that:
\begin{equation}
\label{Abelfirst-2}
k(y) =  k =\int_0^y \frac{f(z)}{(y-z)^{1/2}} dz.
\end{equation}
Now, we proceed in a different manner, multiplying both sides
of~(\ref{Abelfirst-2}) by $1/\Gamma(1/2) = 1/\sqrt{\pi}$ in order to obtain
\begin{eqnarray}
\label{Abel2}
\frac{k}{\Gamma(1/2)} &=& \frac{1}{\Gamma(1/2)} \int_0^y \frac{f(z)}{(y-z)^{1/2}} dz\nonumber \\
&=& _{0}{\textrm{I}}_z^{1/2} \left[f(z) \right],
\end{eqnarray}
whose right-hand side we can recognize (better: we define!)  as a semi-integral operator, that is,   a fractional integral operator of order $1/2$, represented by the symbol ${}_0{\textrm{I}}_z^{1/2}$ we have introduced before.  If this is true, to cancel the ``semi-operator'', we have just to take the ``semi-derivative'' of each side of~(\ref{Abel2}), that is,
\begin{equation}
\label{eq*}
\frac{d^{1/2}}{dz^{1/2}}\left[\frac{k}{\Gamma(\alpha)} \right]= \frac{d^{1/2}}{dz^{1/2}} \bigg\{{}_0{\textrm{I}}_z^{1/2} \left[f(z) \right] \bigg\}= f(z).
\end{equation}
In this way, the problem is reduced to calculate the derivative of order $1/2$ of the left-hand side, i.e., of a constant. This calculation may be performed by using a result due to Lacroix~\cite{Lacroix1819}.
 
Indeed, a derivative of arbitrary order was briefly mentioned (perhaps for the first time as such) in 1819, in the book of Sylvestre Fran\c{c}ois Lacroix (1765--1843), from which an explicit formula for a fractional derivative may be obtained. The central argument may be summarized as follows (using a today's notation). Let, for instance, $y = x^m$; when $n$ is an integer, one has for any arbitrary $m \in \mathbb{N}$:
\begin{eqnarray}
\label{Lacroix}
d^n y = d^n (x^m) &=& m (m-1)\cdots (m-n+1) x^{m-n} dx^n\nonumber \\
&=& \frac{\Gamma(m+1)}{\Gamma(m-n+1)}x^{m-n} dx^n, \quad n\in \mathbb{N},
\end{eqnarray}
where $\mathbb{N}=\{1, 2, 3, ...\}$ is the set of natural numbers. If we put $m=2$ and $n=1$, then we obtain
\begin{equation*}
d^1 (x^2) = 2 x \,dx,
\end{equation*}
as expected. Now, we may consider $m=1$ and $n=1/2$ and try to answer to the original question. We obtain:
\begin{equation}
\label{Lacroix24}
d^{1/2} x = \frac{\Gamma(2)}{\Gamma(3/2)}x^{1/2} dx^{1/2}= \frac{2}{\sqrt{\pi}}\sqrt{x \, dx}.
\end{equation}
Thus, the result obtained by Lacroix may be put in the usual form as:
\begin{equation}
\label{Lacroix1}
\frac{d^{1/2} x}{dx^{1/2}} = \frac{2 \sqrt{x}}{\sqrt{\pi}},
\end{equation}
which, apart from a constant, may be compared with the result reported by Leibniz, in his letter dated September 30, 1695, to L'H\^opital~\cite{Leibniz1695}:
\begin{displayquote}
{\emph{You see from this, Sir, that one can express by an infinite series a quantity such as $d^{1/2} \overline{xy}$, or $d^{1:2} \overline{xy}$, even though this may seem far removed from Geometry, which ordinarily knows only differences with positive integer exponents, or negative ones with respect to sums, and not yet those whose exponents are fractional. It is true that it remains to provide $d^{1:2} x$ for that series; but even this can be explained in some fashion. For let the ordinates $x$ be in geometric progression, so that, taking a constant $d\beta$, one has $dx = x d\beta : a$, or (taking $a$ as unity) $dx = x d\beta$. Then $ddx$ will be $x \overline{d\beta}^2$, and $dx^3$ will be $= x \overline{d\beta}^3$, etc., and $d^e x = x \cdot \overline{d\beta}^e$. And by this device, the potential exponent, and replacing $dx : x$ with $d\beta$, we have $d^e x = \overline{dx : x}^e \cdot x$. Thus it follows that $d^{1:2} x$ will be equal to $x \cdot \sqrt[2]{dx : x}$. It seems likely that one day useful consequences will be drawn from these seemingly useless paradoxes. You are among those who can go furthest in discoveries, and I shall soon be obliged \textit{to hand over the torch to others}. I wish I had much to share, for this verse: \textit{your knowledge is nothing unless someone else knows that you know} is most true, in that thoughts which were of little value in themselves can give rise to far more beautiful ones.}}
\end{displayquote}
Thus, the result presented by Leibniz for the question raised by L'H\^opital may be interpreted as stating that for the fractional exponent $e$, one has:
\begin{equation}
d^e x = x^{1-e} dx^e,
\end{equation}
which could be compared with Eq.~(\ref{Lacroix}), because, for $m=1$ and $n=e$, it could be rewritten as:

\begin{equation}
d^e x =\frac{\Gamma(2)}{\Gamma(2-e)} x^{1-e} dx^e.
\end{equation}
In this framework, if we consider consider the case $n=e=1/2$ and $m=0$, that is, $y=x^0 =1$, from Eq.~(\ref{Lacroix}), we obtain the ``surprising'' result that the fractional derivative of a constant is not zero, but instead a function of $x$, namely:
\begin{equation}
\frac{d^{1/2}}{dx^{1/2}} 1 = \frac{1}{\sqrt{\pi x}}.
\end{equation}

Coming back to our problem, by using the results above in Eq.~(\ref{eq*}), we obtain:
\begin{equation}
f(z) = \frac{k}{\pi \sqrt{z}}.
\end{equation}
This permits us to deduce that the answer to the tautochrone problem is formally:
\begin{equation*}
\frac{ds}{dy} = \frac{\sqrt{2g}}{\pi} \frac{T(y_0)}{\sqrt{y}}.
\end{equation*}
Integrating the above equation, using the condition $s(0)=0$, we obtain:
\begin{equation}
\label{sdey}
s(y) = 2\sqrt{2 g \left[\frac{T(y_0)}{\pi} \right]^2} y^{1/2},
\end{equation}
which is the cycloid sketched in Fig.~\ref{Fig:Abel}, as shown by Abel in connection with this problem of finding the path followed by the particle in order the time of sliding down be independent of the initial height~\cite{Abel1823}.

In summary, to solve the tautochrone problem, in the pioneering Abel's perspective,  we have to determine $f(z)$, a task that can be accomplished by computing the fractional derivative of a constant $k$. In the paper of Abel, it is not clear that if he used the Lacroix or its own results. Anyway, in the paper of 1826, Abel provided the solution for the integral equation~\cite{Abel1826}:
\begin{equation}
\label{Abelsecond}
k(y) = \int_c^y \frac{f(z)}{(y-z)^{\alpha}}dz, \quad y>c, \quad 0< \alpha <1.
\end{equation}
In practice, Abel used the operators that nowadays are ascribed to Riemann and Liouville, preceding them by at least one decade~\cite{Mainardi2007,Mainardi2007b}. To complete the analysis, let us formalize a little bit more the previous results, considering the \textit{Abel fractional integral equation of first kind}, defined as:
\begin{equation}
\label{Abel-1M}
{}_0{\textrm{I}}_t^{\alpha} g(t) = \frac{1}{\Gamma(\alpha)}\int_0^t \frac{g (\tau)}{(t-\tau)^{1-\alpha}} d\tau = f(t), \quad 0 < \alpha < 1.
\end{equation}
Equation~(\ref{Abelfirst}) is a particular case of Eq.~(\ref{Abel-1M}), for $\alpha = 1/2$.
As we have seen above, Eq.~(\ref{Abel-1M}) may be solved in terms of a fractional derivative in the form:
\begin{equation}
\label{DI}
{}_0{\textrm{D}}_t^{\alpha} \left[{}_0{\textrm{I}}_t^{\alpha} g(t)\right] = g(t) = {}_0{\textrm{D}}_t^{\alpha} f(t),
\end{equation}
where we have also introduced the symbol ${}_0{\textrm{D}}_t^{\alpha}$ now to denote a derivative of arbitrary order $\alpha$ in the variable $t$. A notation of this kind was introduced by the mathematician Harold Thayer Davis (1892--1974)~\cite{Davis} and is useful to highlight the meaning of the left subscript, which is not obvious in a derivative~\cite{Mainardi-Book}.

In Eq.~(\ref{DI}), we have supposed that ${}_0{\textrm{D}}_t^{\alpha}$ exists (even if we have not discovered its form until now!) and is the left-inverse of
the operator ${}_0{\textrm{I}}_t^{\alpha}$, that is,
\begin{equation}
\label{Identity}
{}_0{\textrm{D}}_t^{\alpha} {}_0{\textrm{I}}_t^{\alpha} = \mathds{1},
\end{equation}
where $\mathds{1}$ is the identity operator. A solution for~(\ref{Abel-1M}) may be now searched by using the Laplace transform and the convolution theorem~\cite{Mainardi-Book,Evangelista2023}.

\section{The Fractional Operators of Riemmann-Liouville and Caputo}

The expression \emph{fractional calculus} we have been using until now --- although something of a misnomer --- refers to the theory concerning integrals and derivatives of arbitrary order. In this context, it may be viewed as a natural extension of the classical concepts of differentiation and integration to non-integer (fractional) orders. A conventional way to introduce the basic ideas of fractional calculus is by posing a question that has become emblematic: what meaning can we assign to the expression
\begin{equation}
\label{Equation1Frac}
\frac{d^{1/2} f(x) }{dx^{1/2}}?
\end{equation}
This very question was famously raised by Guillaume Fran{\c c}ois Antoine, marquis de L'H\^opital (1661--1704), in a query addressed to Gottfried Wilhelm Leibniz (1646--1716), the founder of differential calculus, whose answer we quoted after Eq.~(\ref{Lacroix24}) in discussing the Lacroix derivative~\cite{Leibniz1695}.

The foundational questions that marked the birth of fractional calculus --- once seen as paradoxical or speculative --- have now been thoroughly examined and answered in full detail. These inquiries reveal profound conceptual challenges. For instance, if the first derivative of a function represents its instantaneous rate of change or the slope at a point, what geometric meaning can we assign to a derivative of order $1/2$? And if we persist in exploring the notion of half-order differentiation, we are naturally led to ask: what operator, when applied twice, yields the first derivative? Such questions, once open-ended, have been rigorously addressed, and the field of fractional calculus has since matured into a well-established and robust branch of mathematical analysis~\cite{Oldham1974,Podlubny,Samko}.

Today, the Riemann-Liouville fractional integral or Riemann-Liouville integral of arbitrary order $\alpha$ is defined as:
\begin{equation}
\label{RLI}
{}_c{\textrm{D}}_x^{-\alpha} f(x) = \frac{1}{\Gamma(\alpha)} \int_c^x \frac{f(t)}{(x-t)^{1-\alpha}} dt, \quad \quad \Re(\alpha) > 0,
\end{equation}
where $\Re(\dots)$  denotes the real part of $\dots$. Then, it is clear that the fractional derivative of arbitrary order $\alpha$ may be interpreted as the derivative of positive integer order of a derivative of negative non-integer order, which is an integral of arbitrary non-integer order. Accordingly, the operator
\begin{eqnarray}
\label{RLIP}
{}_c{\textrm{D}}_x^{\alpha}f(x) = \frac{1}{\Gamma(k-\alpha)}\frac{d^k}{dx^k}\left[\int_c^x \frac{ f(t)}{(x-t)^{\alpha + 1 -k}} dt\right], \quad k \in \mathbb{N},
\end{eqnarray}
for $\alpha = k - p$ is the definition of the Riemann-Liouville fractional derivative to be used hereafter.

A concise, alternative way, to face the problem of using the definition of an integral of arbitrary order to obtain a definition of a derivative of arbitrary order may be sketched as follows. We start with the Riemann-Liouville integral defined from the Cauchy integral formula, as before:
\begin{equation}
\label{A1}
{}_c{\textrm{I}}_x^{\alpha} f(x) = \frac{1}{\Gamma(\alpha)}\int_c^x (x-t)^{\alpha-1} f(t) dt.
\end{equation}
Notice that $\alpha$ may be complex, with the real part strictly positive, that is,  $\Re(\alpha)>0$. We assume for a moment that $\alpha \in \mathbb{R}$. The definition~(\ref{A1}) has the following properties:
\begin{eqnarray*}
{}_c{\textrm{I}}_x^{\alpha} \left[{}_c{\textrm{I}}_x^{\beta}f(x)\right] = {}_c{\textrm{I}}_x^{\alpha+\beta} f(x) \quad {\textrm{and}}\quad \frac{d}{dx} \left[{}_c{\textrm{I}}_x^{\alpha +1}f(x)\right] = {}_c{\textrm{I}}_x^{\alpha}f(x).
\end{eqnarray*}
One is then tempted to assume that the fractional derivative could be easily defined as
\begin{equation*}
{}_c{\textrm{D}}_x^{\alpha}f(x) = {}_c{\textrm{I}}_x^{-\alpha}f(x).
\end{equation*}
One of the problems with this identification is that the gamma function is not defined for zero or negative integers. Anyway, we may follow a procedure inspired in this idea with some caution. We notice that for $n \in \mathbb{N}$, we surely write:
\begin{equation*}
\frac{d^n}{dx^n} \left[{}_c{\textrm{I}}_x^{n}f(x) \right]= f(x),
\end{equation*}
that is, taking $n$ times the derivative of a function $f(x)$ after integration it $n$ times is equivalent to the identity operator, as stated in~(\ref{Identity}). With this result in mind, we expect to find a fractional derivative operator such that
\begin{equation*}
{}_c{\textrm{D}}_x^{\alpha} \left[{}_c{\textrm{I}}_x^{\alpha}f(x)\right] = f(x).
\end{equation*}
An operator that can be constructed is
\begin{equation*}
{}_c{\textrm{D}}_x^{\alpha}f(x) = \frac{d^{\lceil{\alpha}\rceil}}{dx^{\lceil{\alpha}\rceil}} \left[{}_c{\textrm{I}}_x^{{\lceil{\alpha}\rceil}-\alpha}f(x)\right],
\end{equation*}
with ${\lceil{\alpha}\rceil}$ being the ceiling function, which gives the smallest integer greater or equal to $\alpha$ (that is, one has to consider the next integer). This permits us to arrive at the operator defined in~(\ref{RLIP}) for $\alpha = k - p$, such that $0< p=k-\alpha <1$, and ${}_c{\textrm{D}}_x^{k} = d^k/dx^k$ being the usual derivative of integer order. It is, indeed, the \textit{left} Riemann-Liouville fractional derivative~\cite{Yang-Gao-Ju}.

Along the same lines, the Caputo operator may be defined as follows~\cite{Caputo1967,Caputo1969}:

\begin{equation}
\label{Caputo-2}
{}^{\textrm{C}}_c{\textrm{D}}_x^{\alpha} f(x) = \frac{1}{\Gamma(m-\alpha)} \int_c^x \frac{f^{(m)}(t)}{(x-t)^{\alpha +1 -m}}dt, \quad \quad m-1< \alpha < m,
\end{equation}
which, for $m=1$, becomes:
\begin{equation}
\label{Caputo-3}
{}^{\textrm{C}}_c{\textrm{D}}_x^{\alpha} f(x) = \frac{1}{\Gamma(1-\alpha)} \int_c^x \frac{f'(t)}{(x-t)^{\alpha }}dt, \quad \quad 0< \alpha < 1.
\end{equation}
Apart the constant coefficient represented by $1/\Gamma(1-\alpha)$, the operator defined in Eq.~(\ref{Caputo-3}) coincides with the definition used by Abel --- the one we have rewritten in Eq.~(\ref{preCaputo}). To conclude the abbreviated presentation of these operators, we remark that 
a relation between the Caputo and the Riemann-Liouville derivatives may be established as follows~\cite{Podlubny}:
\begin{equation}
{}^{\textrm{C}}_0{\textrm{D}}_x^{\alpha} f(x) = {}_0{\textrm{D}}_x^{\alpha} f(x) -\sum_{k=0}^{n-1}
\frac{f^{(k)}(0) x^{k-\alpha}}{\Gamma(k-\alpha +1)}, \quad n-1 < \alpha < n,
\end{equation}
which, when $0 <\alpha <1$, yields
\begin{eqnarray}
\label{EqC}
{}^{\textrm{C}}_0{\textrm{D}}_x^{\alpha} f(x) &=& {}_0{\textrm{D}}_x^{\alpha} \left[f(x) - f(0)\right] \nonumber \\
&=&  {}_0{\textrm{D}}_x^{\alpha} f(x) - \frac{f(0) x^{-\alpha}}{\Gamma(1-\alpha)}.
\end{eqnarray}
We notice that, from the above equation, Caputo and Riemann-Liouville operators coincide for the functions such that $f(0) =0$. In addition, if we put now $c=-\infty$ in the Riemann-Liouville derivative, then it becomes
\begin{eqnarray}
{}_{-\infty}{\textrm{D}}_x^{\alpha} f(x) &=& \frac{d^m}{dx^m}\left[\frac{1}{\Gamma(m-\alpha)}\int_{-\infty}^x \frac{f(t)}{(x-t)^{\alpha + 1 - m}}dt\right] \nonumber \\
&=& \frac{1}{\Gamma(m-\alpha)} \int_{-\infty}^x
\frac{f^{(m)}(t)}{(x-t)^{\alpha + 1 - m}}dt \nonumber \\
&=&{{}_{-\infty}\!\!}{}^{\textrm{C}}{\textrm{D}}_x^{\alpha} f(x), \quad \quad m-1 < \alpha < m.
\end{eqnarray}
In this limit, both definitions become equal provided that $f(x)$ and its derivatives have a reasonable behavior when $x\to \infty$, that is, $f^{(k)}(-\infty)\to 0$, for
$k = 0, 1, ..., n-1$, with $n = {\lceil{\alpha}\rceil}$.   This is a very important property from the physical point of view because to consider stationary processes is now permitted as, for instance, in the response of the fractional order dynamical systems to a periodic signal, required in the impedance problems~\cite{Evangelista2018}, in the wave propagation in continuous media~\cite{Mainardi-Book,Mainardi-C}, etc.

\section{Concluding Remarks}

In this work, we have  highlighted a convergence of interests between Huygens and Abel ---separated by nearly two centuries, yet united in their engagement with mathematical problems deeply connected to physics. For this analysis, we employed the tautochrone problem --- one of the significant isoperimetric problems in the scientific literature --- as a guiding thread.

We have outlined how the tautochrone was initially studied by Huygens in the 17th century, motivated by physical considerations related to timekeeping, and how its solution involved the cycloid curve. Later, in the 19th century, Niels Henrik Abel, by considering the inverse problem (reconstructing the curve from the descent time as a function of position), arrived at an integral equation --- of the Volterra type --- which can be reinterpreted, in more modern terms, as a particular case of a fractional equation.

This transition from a classical physical-geometric problem to an abstract formulation involving generalized integrals is a compelling example of the evolution of mathematical language and the conceptual tools of mathematical physics. It is of particular interest to observe how Abel, already in the early 19th century, worked with integrals of arbitrary order, anticipating aspects of fractional calculus and laying the groundwork for the consolidation of its foundations --- whose relevance is attested by the interests of our own time. This relevance lies both in the more formal aspects of mathematical analysis and in the numerous applications of this formalism across various scientific domains.

\subsection*{Acknowledgments}
L. R. E. gratefully acknowledges the hospitality of the Department of Molecular Sciences and Nanosystems at Ca' Foscari University of Venice, Italy, and the Department of Materials, Environmental Sciences and Urban Planning (SIMAU) at the Polytechnic University of Marche (UnivPM), Ancona, Italy, where he served as a Visiting Scholar during the preparation of this work.
\\
The research of F. M. has been carried out in the framework 
of the activities of 
the Italian National Group of Mathematical Physics (GNFM-INdAM).
\\
The authors are grateful to the anonymous referees for valuable suggestions which help us to improve the presentation of the results.


\subsection*{Financial Support}
No specific financial support was required for the work. L. R. E. acknowledges the continuous partial financial support from the Brazilian Agency --- CNPq, by means of the Grant n. 30295/2022-9.

\subsection*{Authors' Contributions}

Both authors contributed equally to the final version of this paper. The initial idea was suggested by F.M. The conceptual framework was developed after a fruitful meeting in Pisa in 2024. The first draft was written by L.R.E. Refinements and the first complete version emerged after an extended and productive exchange of correspondence between the two authors. The current and final form of the paper is the result of a joint effort, consolidated during a personal meeting in Rimini --- land of Francesca and Paolo Malatesta, whose tragic love was immortalized by Dante in the Divine Comedy --- on the shores of the Adriatic Sea.

\subsection*{Conflict of Interest}

The authors declare no potential competing interests.

\subsection*{Data Availability Statement} 

Data sharing is not applicable to this article as no new data were created or analyzed in this study.

\subsection*{Mathematics Subject Classification.}   
\begin{itemize}
\item 01 --- History and biography.
\item 01 --- 02 Research exposition pertaining to history and biography
\end{itemize}

\subsection*{Generative Artificial Intelligence (AI)}
This manuscript benefited from the use of Generative AI tools (ChatGPT) for translation assistance and grammar refinement, given that the authors are not native speakers of English.

\subsection*{Authors' Biography}
{\bf Luis Roberto Evangelista} is Professor of Theoretical Physics at the PhD Preogram in Physics at the University of Maring\'a (UEM), Brazil. His research interest lie in complex fluids, complex systems, and history of physics, and include mathematical physics of liquid-crystals, diffusion problems, and adsorption-desorption phenomena. 
\\
{\bf Francesco Mainardi} is Retired Professor  of Mathematical Physics at the University of Bologna, Italy, renowned for his pioneering work in fractional calculus.
His research focuses on fractional differential equations, viscoelasticity, and anomalous diffusion.
He played a key role in promoting the Mittag-Leffler and Wright functions in applied contexts.
Mainardi authored influential books, including Fractional Calculus and Waves in Linear Viscoelasticity.
He remains an active voice in the field, shaping theory and applications across physics and engineering.

\vfill\eject

\end{document}